\documentclass[english,oneside]{smfart}

\newcommand{\opn}[1]{\operatorname{#1}}

\usepackage{amsthm,mathtools,amssymb,textcomp,extarrows,bm,mleftright,graphicx,stmaryrd,scalerel,tikz,enumitem,float,booktabs}
\usetikzlibrary{cd}
\mleftright

\makeatletter
\AtBeginDocument{%
 \let\glb@currsize\relax
}
\makeatother

\usepackage{mathrsfs}

\ExplSyntaxOn
\NewDocumentCommand{\definealphabet}{mmmm}
 {%
  \int_step_inline:nnn { `#3 } { `#4 }
   {
    \cs_new_protected:cpx { #1 \char_generate:nn { ##1 }{ 11 } }
     {
      \exp_not:N #2 { \char_generate:nn { ##1 } { 11 } }
     }
   }
 }
\ExplSyntaxOff
\definealphabet{bb}{\mathbb}{A}{Z}
\definealphabet{cal}{\mathcal}{A}{Z}
\definealphabet{frak}{\mathfrak}{A}{z}
\definealphabet{rm}{\mathrm}{A}{z}
\definealphabet{bf}{\mathbf}{A}{Z}
\definealphabet{scr}{\mathscr}{A}{Z}
\definealphabet{tt}{\mathtt}{A}{Z}
\DeclareMathOperator{\Gal}{Gal}

\newcommand{\lto}{\longrightarrow}

\newcommand{\wtilde}{\widetilde}

\newcommand{\bbrac}[1]{\llbracket #1 \rrbracket}

\newcommand{\hyphen}{\!\opn{-}\!}
\newcommand{\ulpi}{\underline{\pi}}
\usepackage[pdfencoding=auto,psdextra]{hyperref}
\usepackage[nameinlink]{cleveref}
\Crefname{theorem}{Theorem}{Theorem}
\Crefname{conjecture}{Conjecture}{Conjectures}
\Crefname{lemma}{Lemma}{Lemmas}
\Crefname{definition}{Definition}{Definitions}
\Crefname{remark}{Remark}{Remarks}
\Crefname{proposition}{Proposition}{Propositions}
\Crefname{corollary}{Corollary}{Corollaries}
\Crefname{equation}{}{}
\Crefname{item}{}{}
\Crefname{example}{Example}{Examples}
\Crefname{proof}{Proof}{Proofs}
\newlist{thmenum}{enumerate}{1}
\setlist[thmenum]{label=\arabic*., ref=\thetheorem~(\arabic*)}
\crefalias{thmenumi}{theorem}
\newlist{propenum}{enumerate}{1}
\setlist[propenum]{label=\arabic*., ref=\theproposition~(\arabic*)}
\crefalias{propenumi}{proposition}
\newlist{lemenum}{enumerate}{1}
\setlist[lemenum]{label=\arabic*., ref=\thelemma~(\arabic*)}
\crefalias{lemenumi}{lemma}
\usepackage[url=true,backend=biber,style=ext-alphabetic,hyperref=true,giveninits=true]{biblatex}
\makeatletter
\renewcommand*\subjclass[2][2020]{%
  \def\@subjclass{#2}%
  \@ifundefined{subjclassname@#1}{%
    \ClassWarning{\@classname}{Unknown edition (#1) of Mathematics
      Subject Classification; using '2020'.}%
  }{%
    \@xp\let\@xp\subjclassname\csname subjclassname@#1\endcsname
  }%
}
\@namedef{subjclassname@1991}{%
  \textup{1991} Mathematics Subject Classification}
\@namedef{subjclassname@2000}{%
  \textup{2000} Mathematics Subject Classification}
\@namedef{subjclassname@2010}{%
  \textup{2010} Mathematics Subject Classification}
\@namedef{subjclassname@2020}{%
  \textup{2020} Mathematics Subject Classification}
\@xp\let\@xp\subjclassname\csname subjclassname@2020\endcsname
\makeatother
\usepackage{framed}
\newtheorem{theorem}{Theorem}[section]
\newtheorem{example}[theorem]{Example}
\AddToHook{env/example/begin}{\crefalias{theorem}{example}}
\newtheorem{lemma}[theorem]{Lemma}
\AddToHook{env/lemma/begin}{\crefalias{theorem}{lemma}}
\newtheorem{remark}[theorem]{Remark}
\AddToHook{env/remark/begin}{\crefalias{theorem}{remark}}
\newtheorem{proposition}[theorem]{Proposition}
\AddToHook{env/proposition/begin}{\crefalias{theorem}{proposition}}
\newtheorem{definition}[theorem]{Definition}
\AddToHook{env/definition/begin}{\crefalias{theorem}{definition}}

\AddToHook{env/conjecture/begin}{\crefalias{theorem}{conjecture}}

\AddToHook{env/corollary/begin}{\crefalias{theorem}{corollary}}

\AddToHook{env/question/begin}{\crefalias{theorem}{question}}

\AddToHook{env/assumption/begin}{\crefalias{theorem}{assumption}}

\newtheorem{introthm}{Theorem}

\AddToHook{env/theorem/before}{\addvspace{3pt}}
\AddToHook{env/theorem/after}{\addvspace{3pt}}
\AddToHook{env/example/before}{\addvspace{3pt}}
\AddToHook{env/example/after}{\addvspace{3pt}}
\AddToHook{env/definition/before}{\addvspace{3pt}}
\AddToHook{env/definition/after}{\addvspace{3pt}}
\AddToHook{env/proposition/before}{\addvspace{3pt}}
\AddToHook{env/proposition/after}{\addvspace{3pt}}
\AddToHook{env/lemma/before}{\addvspace{3pt}}
\AddToHook{env/lemma/after}{\addvspace{3pt}}
\AddToHook{env/conjecture/before}{\addvspace{3pt}}
\AddToHook{env/conjecture/after}{\addvspace{3pt}}
\AddToHook{env/introthm/before}{\addvspace{3pt}}
\AddToHook{env/introthm/after}{\addvspace{3pt}}
\AddToHook{env/introconj/before}{\addvspace{3pt}}
\AddToHook{env/introconj/after}{\addvspace{3pt}}
\AddToHook{env/remark/before}{\addvspace{3pt}}
\AddToHook{env/remark/after}{\addvspace{3pt}}
\AddToHook{env/corollary/before}{\addvspace{3pt}}
\AddToHook{env/corollary/after}{\addvspace{3pt}}
\Crefformat{enumi}{#2\textup{(#1)}#3}
\usepackage{microtype}
\usepackage{geometry}
\definecolor{labelkey}{rgb}{1,0,0}
\definecolor{refkey}{rgb}{0,0,1}

\title{Descent of $(\varphi,\tau)$-modules in characteristic $p$}

\usepackage{orcidlink}
\author{Yijun Yuan\orcidlink{0000-0001-6571-6980}}
\address{Institute for Theoretical Sciences, Westlake University, No. 600 Dunyu Road, Sandun town, Xihu district, Hangzhou, Zhejiang Province, 310030, China}
\email{941201yuan@gmail.com}
\urladdr{https://yijunyuan.github.io/}

\begin{document}
\frontmatter
\begin{abstract}
	In this article, we study the descent of $(\varphi,\tau)$-modules over perfectoid period rings in characteristic $p$ via Berger and Rozensztajn's theory of super-H\"{o}lder vectors. This is a generalization of their work on $(\varphi,\Gamma)$-modules. As an application, we answer a question of Caruso regarding the connection between $(\varphi,\tau)$-modules and $(\varphi,\Gamma)$-modules without involving Galois representations as intermediaries.
\end{abstract}
\subjclass{11S15, 13J05, 22E35}
\keywords{$(\varphi,\tau)$-modules, super-H\"{o}lder vectors, decompletion, deperfection}
\maketitle
\tableofcontents
\mainmatter
\section{Introduction}
Let $p\geq 3$ be a prime number. In 1990, J. M. Fontaine introduced the concept of étale \((\varphi,\Gamma)\)-modules (cf. \cite{Fontaine2007}), which constitute a class of modules over period rings endowed with both a Frobenius action \(\varphi\) and an action of the group \(\Gamma_K := \mathrm{Gal}(K(\zeta_{p^\infty})/K)\). The category of \((\varphi,\Gamma)\)-modules is equivalent to the category of \(p\)-adic Galois representations of \(\mathbf{Q}_p\) (cf. \cite[Théorème 3.4.3]{Fontaine2007}), and has since become a foundational tool in the study of \(p\)-adic Hodge theory.

As the theoretical framework has developed, it has become increasingly recognized that descending $(\varphi,\Gamma)$-modules from large period rings to smaller ones is a fundamental process. For instance, F. Cherbonnier and P. Colmez established that every $(\varphi,\Gamma)$-module defined over Fontaine's ring $\mathbf{A}_K$ is overconvergent (cf. \cite{cherbonnierRepresentationsPadiquesSurconvergentes1998}). Moreover, the theory of locally analytic vectors offers a more systematic approach to addressing issues of overconvergence and facilitates the process of decompletion.

However, the classical theory of locally analytic vectors only works for $p$-adic representations with rational coefficients, which prevents us from descending the $(\varphi,\Gamma)$-modules with integral or $\opn{mod} p$ coefficients. In the very recent literature, two solutions for resolving this issue have emerged:
\begin{enumerate}
	\item In \cite{bergerDecompletionCyclotomicPerfectoid2022} and  \cite{bergerSuperHolderVectorsField2024}, L. Berger and S. Rozensztajn introduced the notion of super-H\"{o}lder vectors (cf. \Cref{def:47537}), which serves as a $\opn{mod} p$ analogue of the locally analytic vectors. They demonstrated that this framework can be employed to decomplete the field of norms as well as the $(\varphi,\Gamma)$-modules over the perfectoid period ring in characteristic $p$:
	      \begin{theorem}[cf. {\cite[Theorem A]{bergerDecompletionCyclotomicPerfectoid2022}, \cite[Theorem 2.9]{bergerDecompletionCyclotomicPerfectoid2022}}]\label{thm:58509}
		      Let $L/K$ be a totally ramified extension with Galois group $\Gamma$ a $p$-adic Lie group. Let $\wtilde{\bfE}_L\coloneqq \widehat{L}^\flat$ and let $X_K(L)$ be the field of norms of this extension. Then
		      \begin{enumerate}
			      \item One has $\wtilde{\bfE}_L^{\Gamma\hyphen d\hyphen\opn{sh}}=\bigcup_n\varphi^{-n}(X_K(L))$, where $d=\opn{dim}(G)$ and
			            $\wtilde{\bfE}_L^{\Gamma\hyphen d\hyphen\opn{sh}}$ is the set of super-H\"{o}lder vectors in $\wtilde{\bfE}_L$.
			      \item Let $\wtilde{\bfE}_K\coloneqq \widehat{K(\zeta_{p^\infty})}^\flat$ and $\bfE_K\coloneqq X_K(K(\zeta_{p^\infty}))$. Then for any $k,n\in\bfN$, one has $\wtilde{\bfE}_K^{\Gamma^{(k)}\hyphen 1\hyphen\opn{sh},k-n}=\varphi^{-n}(\bfE_K)$, where $\Gamma_K^{(k)}$ is the subgroup of $\Gamma_K$ corresponding to $1+p^k\bfZ_p$ under the cyclotomic character.
		      \end{enumerate}
	      \end{theorem}
	      \begin{theorem}[cf. {\cite[Corollary 3.11]{bergerDecompletionCyclotomicPerfectoid2022}}]\label{thm:45774}
		      Let $D$ be an \'etale $(\varphi,\Gamma)$-module over $\bfE_K$. Then for any $k,n\in\bfN$, one has $\left(\wtilde{\bfE}_K\otimes_{\bfE_K}D\right)^{\Gamma^{(k)}\hyphen 1\hyphen\opn{sh},k-n}=\varphi^{-n}(\bfE_K)\otimes_{\bfE_K}D$.
	      \end{theorem}
	\item In \cite{Porat_2025}, G. Porat developed the theory of mixed-characteristic locally analytic vectors, which encompasses the classical definition and also specializes to super-H\"{o}lder vectors in characteristic $p$. As an application, he reproved \cite[Theorem 2.9]{bergerDecompletionCyclotomicPerfectoid2022} with characteristic $0$ methods.
\end{enumerate}

Parallel to the theory of $(\varphi,\Gamma)$-modules, X. Caruso introduced the theory of $(\varphi,\tau)$-modules (cf. \cite{carusoRepresentationsGaloisiennesPadiques2013}), which makes use of the Kummer extension $K\left(\pi^{1/p^\infty}\right)/K$ instead of the cyclotomic extension $K(\zeta_{p^\infty})/K$. This variant is naturally connected to the Breuil-Kisin modules, which plays an important role in integral $p$-adic Hodge theory. Since the Kummer extension is not Galois over $K$, One has to enlarge it to the period ring of its normal closure to define a reasonable Galois action. More precisely, we have the following definition:
\begin{definition}[cf. {\cite[Definition 6.2.2]{gaoLOCALLYANALYTICVECTORS2021}}]
	Let $K_\infty\coloneqq \bigcup_n K\left(\pi^{1/p^n}\right)$ be the Kummer extension of $K$ and let $K_\rtimes$ be the normal closure of $K_\infty/K$. Let $\wtilde{\bfE}_{\rtimes,K}\coloneqq \widehat{K_\rtimes}^\flat$, $\wtilde{\bfE}_{\tau,K}\coloneqq \widehat{K_\infty}^\flat$ and $\bfE_{\tau,K}\coloneqq X_K(K_\infty)$. For $\bfE_{?}\in\{\bfE_{\tau,K},\wtilde{\bfE}_{\tau,K}\}$, an \textbf{étale $(\varphi,\tau)$-module} over $\left(\bfE_{?},\wtilde{\bfE}_{\rtimes,K}\right)$ consists of
	\begin{enumerate}
		\item An \'etale $\varphi$-module $M$ over $\bfE_{?}$, i.e. a finite dimensional vector space over $\bfE_{\tau,K}$ equipped with a semilinear Frobenius action $\varphi\colon M\lto M$ such that $\varphi^*(M)\cong M$;
		\item A semilinear action of $\Gamma_{\rtimes,K}\coloneqq \opn{Gal}\left(K_\rtimes/K\right)$ on $\wtilde{\bfE}_{\rtimes,K}\otimes_{\bfE_{?}}M$ that commutes with $\varphi$ on $\wtilde{\bfE}_{\rtimes,K}\otimes_{\bfE_{?}}M$, satisfying $M\subseteq \left(\wtilde{\bfE}_{\rtimes,K}\otimes_{\bfE_{?}}M\right)^{\opn{Gal}(K_\rtimes/K_\infty)}$.
	\end{enumerate}
\end{definition}

In spirit of \Cref{thm:58509} and \Cref{thm:45774}, it is natural to ask whether one can decomplete and deperfect the field of norms $\bfE_{\tau,K}$ and the $(\varphi,\tau)$-modules over $\left(\wtilde{\bfE}_{\tau,K},\wtilde{\bfE}_{\rtimes,K}\right)$ by using super-H\"{o}lder vectors? Since there is no well-defined ``$\tau$-action'' on the underlying $\varphi$-module of a $(\varphi,\tau)$-module, a direct search-and-replace rewrite of the arguments in \cite{bergerDecompletionCyclotomicPerfectoid2022} does not work. On the other hand, during their proof of the overconvergence of $(\varphi,\tau)$-modules (over the mixed characteristic period rings), H. Gao and L. Poyeton show that taking $\Gal(K(\zeta_{p^\infty})/K)$-invariant $\Gal(K_\rtimes/K(\zeta_{p^\infty}))$-(locally) analytic vectors on certain period rings will produce the desired decompletion. This motivates us to consider the $\Gal(K(\zeta_{p^\infty})/K)$-invariant $\Gal(K_\rtimes/K(\zeta_{p^\infty}))$-super-H\"{o}lder vectors in our characteristic $p$ setting. To this end, we prove the following theorems:
\begin{introthm}[cf. {\Cref{thm:56067}}]
	One has
	\begin{enumerate}
		\item $$\wtilde{\bfE}_{\rtimes,K}^{\tau\hyphen\opn{sh};\gamma=1}=\wtilde{\bfE}_{\rtimes,K}^{\tau\hyphen\opn{la};\gamma=1}=\varphi^{-\infty}(\bfE_{\tau,K}),$$
		      where $\varphi^{-\infty}(\bfE_{\tau,K})\coloneqq \bigcup_{n\geq 0}\varphi^{-n}(\bfE_{\tau,K})$ is the perfection of $\bfE_{\tau,K}$.
		\item For any $n,k\in\bfN$, one has $\wtilde{\bfE}_{\rtimes,K}^{\tau_k\hyphen\opn{sh},k-n+c_p;\gamma=1}=\varphi^{-n}(\bfE_{\tau,K})$, where $c_p\coloneqq \log_p\frac{p}{p-1}$.
	\end{enumerate}
\end{introthm}
\begin{introthm}[cf. {\Cref{thm:7812}}]\label{thm:B}
	Let $M$ be an \'etale $(\varphi,\tau)$-module over $(\bfE_{\tau,K},\wtilde{\bfE}_{\rtimes,K})$ and $k\in\bfN$. Then
	$$\left(\wtilde{\bfE}_{\rtimes,K}\otimes_{\bfE_{\tau,K}}M\right)^{\tau\hyphen\opn{sh};\gamma=1}=\left(\wtilde{\bfE}_{\rtimes,K}\otimes_{\bfE_{\tau,K}}M\right)^{\tau\hyphen\opn{la};\gamma=1}=\varphi^{-\infty}(\bfE_{\tau,K})\otimes_{\bfE_{\tau,K}}M.$$
	More precisely, for any $n\in\bfN$, one has
	\begin{equation*}
		\left(\wtilde{\bfE}_{\rtimes,K}\otimes_{\bfE_{\tau,K}}M\right)^{\tau_k\hyphen\opn{sh},k-n+c_p;\gamma=1}=\varphi^n(\bfE_{\tau,K})\otimes_{\bfE_{\tau,K}}M.
	\end{equation*}
\end{introthm}

\begin{remark}\leavevmode
	\begin{enumerate}
		\item For an \'etale $(\varphi,\tau)$-module $M$ over $(\bfE_{\tau,K},\wtilde{\bfE}_{\rtimes,K})$, $\wtilde{\bfE}_{\tau,K}\otimes_{\bfE_{\tau,K}}M$ is the \'etale $(\varphi,\tau)$-module over $(\wtilde{\bfE}_{\tau,K},\wtilde{\bfE}_{\rtimes,K})$ that associated to the same Galois representation as $M$. If we write $\wtilde{\bfE}_{\rtimes,K}\otimes_{\bfE_{\tau,K}}\cong\wtilde{\bfE}_{\rtimes,K}\otimes_{\wtilde{\bfE}_{\tau,K}}\left(\wtilde{\bfE}_{\tau,K}\otimes_{\bfE_{\tau,K}}M\right)$, then \Cref{thm:B} can be viewed as a decompletion result for \'etale $(\varphi,\tau)$-modules over perfectoid period rings $(\wtilde{\bfE}_{\tau,K},\wtilde{\bfE}_{\rtimes,K})$.

		\item The notation ``$\tau\hyphen\opn{la}$'' in the above theorems refers to Porat's locally analytic vectors, which coincide with the super-H\"{o}lder vectors in our characteristic $p$ setting (cf. \cite[Example 2.14 (2)]{Porat_2025}).
	\end{enumerate}
\end{remark}

\subsection*{New difficulties and new methods in the Kummer case}
Compared to the case of $(\varphi,\Gamma)$-modules, the proofs of the above theorems involve several new ingredients:
\begin{enumerate}
	\item In \cite{bergerSuperHolderVectorsField2024}, to apply the Ax-Sen-Tate theorem, the authors required the extension to be $p$-adic Lie extension, which allows them to simplify the description of the image of the embedding of the field of norms into the perfectoid field via Sen's theorem (cf. \cite{senRamificationPadicLie1972}). However, since $K_\infty/K$ is not Galois, the image of the embedding have to be described in a more low-level way: analysing the higher ramification of the extension $K_\infty/K$ (cf. \Cref{prop:7578}).
	\item We avoid the use of normalized Tate trace, which does not behave well in the Kummer setting, during the deperfection. Instead, we observe that the decompletion result (i.e. \Cref{it:59278}) is already strong enough to supersede the role of the Tate trace in \cite{bergerDecompletionCyclotomicPerfectoid2022}.
	\item In \cite{bergerDecompletionCyclotomicPerfectoid2022}, the authors chose a $\Gamma_K$-stable $\bfE_K^+$-lattice in a $(\varphi,\Gamma)$-module $D$ to induced a $\Gamma_K$-isometric valuation on $D$. However, we do not have such lattice in the underlying $\varphi$-module of a $(\varphi,\tau)$-module $M$. This difficulty is resolved by considering two different valuations on $\wtilde{\bfE}_{\rtimes,K}\otimes_{\bfE_{\tau,K}}M$ (cf. \Cref{sec:37644}): one is induced from a lattice of the underlying $\varphi$-module, while the other one is $\Gamma_{\tau,K}$-isometric. The interplay between these two valuations, whose difference can be controlled, allows us to adapt the arguments in \cite{bergerDecompletionCyclotomicPerfectoid2022} to our setting.
\end{enumerate}

\subsection*{A question of Caruso}
In \cite[Section 4]{carusoRepresentationsGaloisiennesPadiques2013}, Caruso raised the question of whether one can connect $(\varphi,\tau)$-modules to $(\varphi,\Gamma)$-modules without involving Galois representations as intermediaries. There are two approaches to this problem in the existing literature:
\begin{enumerate}
	\item In L. Poyeton's PhD thesis (cf. \cite[Théorème I]{poyetonExtensionsLiePadiques2019}), he used the theory of pro-analytic vectors to recover $(\varphi,\Gamma)$-modules from the extension of scalars of $(\varphi,\tau)$-modules over the Robba ring to $\wtilde{\bfB}_{\opn{rig},L}^\dagger\coloneqq \left(\wtilde{\bfE}_{\opn{rig}}^\dagger\right)^{\scrG_{K_\rtimes}}$.
	\item In \cite[Theorem B]{yuan2025multivariableperiodringspadic}, the author replaced $\wtilde{\bfB}_{\opn{rig},L}^\dagger$ in Poyeton's work with the multivariable period rings of $K_\rtimes/K$, which is the $p$-Cohen ring of the completed composition of $\bfE_K$ and $\bfE_{\tau,K}$, to establish a similar equivalence. This avoids the subtle calculation of pro-(locally-) analytic vectors, at the cost of less natural structure on the period rings.
\end{enumerate}
With \cite[Corollary 3.11]{bergerDecompletionCyclotomicPerfectoid2022} and \Cref{thm:B} at hand, we can now directly relate $(\varphi,\tau)$-modules and $(\varphi,\Gamma)$-modules over their corresponding field of norms in characteristic $p$:
\begin{introthm}
	Let $V$ be a $\bfF_p$-representation of $\scrG_K$. Let $\bfD(V)$ (resp. $\bfD_\tau(V)$) be the associated \'etale $(\varphi,\Gamma)$-module over $\bfE_K$ (resp. \'etale $(\varphi,\tau)$-module over $\left(\bfE_{\tau,K},\wtilde{\bfE}_{\rtimes,K}\right)$). Then one has
	$$\left(\wtilde{\bfE}_{\rtimes,K}\otimes_{\bfE_K}\bfD(V)\right)^{\tau\hyphen\opn{sh};\gamma=1}=\left(\wtilde{\bfE}_{\rtimes,K}\otimes_{\bfE_K}\bfD(V)\right)^{\tau\hyphen\opn{la};\gamma=1}\cong\varphi^{-\infty}(\bfE_{\tau,K})\otimes_{\bfE_{\tau,K}}\bfD_\tau(V)$$
	and
	$$\left(\wtilde{\bfE}_{\rtimes,K}\otimes_{\bfE_{\tau,K}}\bfD_\tau(V)\right)^{\gamma\hyphen\opn{sh};\tau=1}=\left(\wtilde{\bfE}_{\rtimes,K}\otimes_{\bfE_{\tau,K}}\bfD_\tau(V)\right)^{\gamma\hyphen\opn{la};\tau=1}\cong\varphi^{-\infty}(\bfE_{K})\otimes_{\bfE_{K}}\bfD(V).$$
	More precisely,
	$$\left(\wtilde{\bfE}_{\rtimes,K}\otimes_{\bfE_K}\bfD(V)\right)^{\tau\hyphen\opn{sh},c_p;\gamma=1}\cong\bfD_\tau(V),\ \left(\wtilde{\bfE}_{\rtimes,K}\otimes_{\bfE_K}\bfD(V)\right)^{\gamma\hyphen\opn{sh},0;\tau=1}\cong\bfD(V).$$
\end{introthm}
\begin{proof}
	The result follows from the following standard isomorphism (see, for instance, \cite[ Proposition6.6]{yuan2025multivariableperiodringspadic}):
	$$\wtilde{\bfE}_{\rtimes,K}\otimes_{\bfE_{\tau,K}}\bfD_\tau(V)\cong\left(\wtilde{\bfE}\otimes_{\bfF_p}V\right)^{\scrG_{K_\rtimes}}\cong \wtilde{\bfE}_{\rtimes,K}\otimes_{\bfE_{K}}\bfD(V).$$
\end{proof}

\subsection*{Conventions}
Since this article is heavily inspired by \cite{bergerDecompletionCyclotomicPerfectoid2022}, many arguments are direct adaptations of those in op. cit. to the context of $(\varphi,\tau)$-modules. To make this adaptation process more transparent, we denote certain references with a dagger symbol ($\dagger$) to signify that they were originally formulated for $(\varphi,\Gamma)$-modules, but have been adapted for application to $(\varphi,\tau)$-modules.

\subsection*{Acknowledgements}
This work is partially derived from the author's PhD thesis at Yau Mathematical Sciences Center, Tsinghua University, under the supervision of Lei Fu. The author would like to express his gratitude to Lei Fu for his support throughout the research process. The author also thanks Laurent Berger for his suggestions on considering super-H\"{o}lder vectors. The work is partially written during the author's visit to Renmin University of China, and he thanks Shanwen Wang for hosting and helpful discussions.

\section{Preliminaries on $p$-adic Kummer extension}
\subsection{Extensions and Galois groups}
\begin{definition}
	Let $K$ be a finite extension of $\bfQ_p$, with ring of integer $\calO_K$, uniformizer $\pi$, and residue field $\kappa$.
	\begin{enumerate}
		\item Let $K_{p^\infty}\coloneqq \bigcup_n K(\zeta_{p^n})$ and $\Gamma_K\coloneqq\opn{Gal}(K_{p^\infty}/K)$.
		\item Let $K_n\coloneqq K(\pi^{1/p^n})$, $K_\infty\coloneqq \bigcup_n K_n$ and $\Gamma_{\tau,K}\coloneqq \opn{Gal}(K_\infty K_{p^\infty}/K_\infty)$.
		\item Let $K_{\rtimes}\coloneqq K_\infty K_{p^\infty}$ and $\Gamma_{\rtimes,K}\coloneqq \opn{Gal}(K_\rtimes/K)$.
		\item Let $\wtilde{\bfE}\coloneqq \bfC_p^\flat$, $\widetilde{\bfE}_K\coloneqq \widehat{K}_{p^\infty}^\flat$, $\widetilde{\bfE}_{\tau,K}\coloneqq \widehat{K}_\infty^\flat$, $\widetilde{\bfE}_{\rtimes,K}\coloneqq \widehat{K}_{\rtimes}^\flat$. We add a superscript $^+$ to denote their valuation rings.
		\item Let $\varepsilon\coloneqq (1,\zeta_p,\cdots,\zeta_{p^n},\cdots)$ and $\ulpi\coloneqq (\pi,\pi^{1/p},\cdots,\pi^{1/p^n},\cdots)$  be elements in $\wtilde{\bfE}^+$.
	\end{enumerate}
\end{definition}
We collect some basic properties of the above extensions and Galois groups in the following lemmas, all of which are well-known to experts.
\begin{lemma}[cf. {\cite[Lemma 5.1.2]{liuLatticesSemistableRepresentations2008}}]\label{lem:34973}\leavevmode
	\begin{lemenum}
		\item\label{it:1257} $K_\infty\cap K_{p^\infty}=K$;
		\item\label{it:45121} $\Gamma_{\rtimes,K}\cong \Gamma_{\tau,K}\rtimes \Gamma_K$, $\opn{Gal}(K_{\rtimes}/K_\infty)\cong\Gamma_K$;
		\item One has isomorphisms
		$$c\colon\Gamma_{\tau,K}\xlongrightarrow{\cong}\opn{Gal}\left(\widetilde{\bfE}_K(\ulpi)\middle/\widetilde{\bfE}_K\right)=\left\{f\colon \ulpi\mapsto \ulpi\cdot \varepsilon^a\middle\vert a\in \bfZ_p\right\}\xlongrightarrow{\cong}\bfZ_p,$$
		where we can take the topological generator $\tau$ of $\Gamma_{\tau,K}$ such that $c(\tau)=1$.
		\item Under the isomorphism $\Gamma_{\bfQ_p}\cong\bfZ_p^\times$ via the cyclotomic character, one has $\Gamma_K\cong \mu\times (1+p^{n(K)}\bfZ_p)$, where $\mu$ is a subgroup of $(\bfZ/p\bfZ)^\times$ and $n(K)\in\bfN$ depends on $K$. We denote by $\Gamma_K^{(n)}$ the subgroup corresponding to $1+p^n\bfZ_p$ for any $n\geq n(K)$.
	\end{lemenum}
\end{lemma}

\begin{lemma}\label{lem:23768}
	For any $n\in\bfN$, let $\Gamma_{\tau,K}^{(n)}\coloneqq c^{-1}\left(p^n\bfZ_p\right)$. Then $\Gamma_{\tau,K}^{(n)}=\opn{Gal}(K_{\rtimes}/K_{p^\infty}K_n)$.
\end{lemma}
\begin{proof}
	The assertion is trivial when $n=0$, therefore we suppose $n\geq 1$. One has trivial inclusion $K_nK_{p^\infty}\subseteq \rmH^0\left(\Gamma_{\tau,K}^{(n)},K_\rtimes\right)$.

	When $i=1$, we have
	$$[K_1K_{p^\infty}\colon K_{p^\infty}]\mid \left[\rmH^0\left(\Gamma_{\tau,K}^{(1)},K_\rtimes\right)\colon K_{p^\infty}\right]=\left[\rmH^0\left(\Gamma_{\tau,K}^{(1)},K_\rtimes\right)\colon \rmH^0\left(\Gamma_{\tau,K}^{(0)},K_\rtimes\right)\right]=p.$$
	If $[K_1K_{p^\infty}\colon K_{p^\infty}]=1$, then $\pi^{1/p}\in K_{p^\infty}\cap K_\infty$, which is impossible by \Cref{it:1257}. Therefore $[K_1K_{p^\infty}\colon K_{p^\infty}]=p$ and consequently $K_nK_{p^\infty}= \rmH^0\left(\Gamma_{\tau,K}^{(n)},K_\rtimes\right)$.

	For $n\geq 2$, notice that
	\begin{align*}
		[K_nK_{p^\infty}\colon K_{p^\infty}]= & [K_nK_{p^\infty}\colon K_{n-1}K_{p^\infty}]\cdot [K_{n-1}K_{p^\infty}\colon K_{p^\infty}]                              \\
		=                                     & \left[(K_{n-1})_1(K_{n-1})_{p^\infty}\colon (K_{n-1})_{p^\infty}\right]\cdot [K_{n-1}K_{p^\infty}\colon K_{p^\infty}],
	\end{align*}
	the result follows from induction.
\end{proof}

\subsection{Field of norms and perfectoid fields}
\begin{definition}[{cf. \cite[Section 1]{wintenbergerCorpsNormesCertaines1983}}]
	Let $L$ be an algebraic extension of $K$.
	\begin{enumerate}
		\item We say that $L/K$ is arithmetic profinite (APF), if for any $u\in\bfR_{\geq -1}$, $\scrG_K^u\scrG_L$ is an open subgroup of $\scrG_K$, where $\scrG_K^u$ is the upper numbering ramification filtration of $\scrG_K$.
		\item We say that $L/K$ is strictly APF (sAPF), if
		      $$\liminf_{u\to\infty}\frac{\psi_{L/K}(u)}{\left(\scrG_K^0\colon \scrG_L^0\scrG_K^u\right)}>0,$$
		      where $\psi_{L/K}$ is the Herbrand function of $L/K$.
	\end{enumerate}
\end{definition}

\begin{example}
	The extensions $K_\infty/K$, $K_{p^\infty}/K$ and $K_\rtimes/K$ are all sAPF extensions.
\end{example}

\begin{theorem}[{cf. \cite[Section 2.1]{wintenbergerCorpsNormesCertaines1983}}]
	Let $L/K$ be an sAPF extension. Call
	$$X_K(L)\coloneqq \{0\}\cup\varprojlim_{M\in\scrE(L/K)}M^\times$$
	the \textbf{field of norms} of $L/K$,
	where $\scrE(L/K)$ is the set of finite extensions of $K$ contained in $L$, and the transition maps are given by the norm maps.
	\begin{enumerate}
		\item The addition
		      $$\left(\alpha_M\right)_{M\in\scrE(L/K)}+\left(\beta_M\right)_{M\in\scrE(L/K)}\coloneqq \left(\varprojlim_{N\in\scrE(L/M)}\scrN_{N/M}(\alpha_N+\beta_N)\right)_{M\in\scrE(L/K)}$$
		      and component-wise multiplication make $X_K(L)$ into a field of characteristic $p$.
		\item Let $K_0$ be the maximal unramified subextension of $K$ in $L$. The field $X_K(L)$ is a complete discretely valued field with respect to the valuation
		      $$v\colon (\alpha_M)_{M\in\scrE(L/K)}\longmapsto v_p(x_{K_0})\cdot e_K,$$
		      where $e_K$ is the absolute ramification index of $K$.
		\item The absolute Galois group of $X_K(L)$ is canonically isomorphic to $\scrG_L$.
	\end{enumerate}
\end{theorem}

\subsubsection*{Embed $X_K(L)$ into $\wtilde{\bfE}_L$}
By a result of Gabber-Ramero (cf. \cite[Proposition 6.6.6]{gabberAlmostRingTheory2003}), if the extension $L/K$ is sAPF, then the completion of $L$ is perfectoid. It is very natural to compare the field of norms $X_K(L)$ with the tilt $\wtilde{\bfE}_L\coloneqq \widehat{L}^\flat$ of $\widehat{L}$. The following theorem answers this question:
\begin{theorem}[cf. {\cite[Proposition 4.2.1, Corollaire 4.3.4]{wintenbergerCorpsNormesCertaines1983}}]
	Let $L/K$ be an sAPF extension, $K_1$ be the maximal tamely ramified subextension of $K$ in $L$. For any $n\in\bfN_{\geq 1}$, let
	$$\scrE_n\coloneqq\{M\in\scrE(L/K_1)\colon v_p([M\colon K_1])\geq n\}.$$
	\begin{enumerate}
		\item For any element $(\alpha_E)_{E\in\scrE(L/K)}$ in $X_K(L)$, the net
		      $$(\scrE_n,\subseteq)\lto\widehat{L},\ E\longmapsto\alpha_E^{[E\colon K_1]\cdot p^{-n}}$$
		      converges to an element $x_n$ in $\widehat{L}$.
		\item The map $(\alpha_E)_{E\in\scrE(L/K)}\longmapsto (x_n)_{n\geq 1}$ gives rise to an continuous embedding of $X_K(L)$ into $\wtilde{\bfE}_L$.
		\item The perfection of the image of $X_K(L)$ in $\wtilde{\bfE}_L$ is dense in $\wtilde{\bfE}_L$.
	\end{enumerate}
\end{theorem}
\begin{definition}
	We identify $X_K(L)$ with its image in $\wtilde{\bfE}_L$ via the above embedding. For the extension $K_\infty/K$ (resp. $K_{p^\infty}/K$), we denote $\bfE_{\tau,K}\coloneqq \bfE_{K_\infty/K}$ (resp. $\bfE_K\coloneqq \bfE_{K_{p^\infty}/K}$).
\end{definition}

The embedding of \( X_K(L) \) into \(\widetilde{\mathbf{E}}_L\) admits a more explicit characterization. To establish this, we require the subsequent lemma:
\begin{lemma}[cf. {\cite[Remarque 1.4.2]{wintenbergerCorpsNormesCertaines1983}}]\leavevmode
	\begin{enumerate}
		\item The set of ramification jumps of $L/K$ is a countable unbounded subset $\{b_1<b_2<\cdots<b_n<\cdots\}$ of $\bfR$.
		\item Call $F_n\coloneqq \overline{\bfQ}_p^{{\scrG}_K^{b_n}\scrG_L}$ the $n$-th elementary extension of $L/K$. One has
		      \begin{enumerate}
			      \item $L=\bigcup_{n\geq 1}F_n$;
			      \item $F_1/K$ is the maximal tamely ramified subextension of $L/K$;
			      \item for any $n\geq 1$, the extension $F_{n+1}/F_n$ is a non-trivial finite $p$-extension.
		      \end{enumerate}
	\end{enumerate}
\end{lemma}
\begin{theorem}\label{thm:10788}
	Let $\bfE_{L/K}^+$ (resp. $\wtilde{\bfE}_L^+$) be the valuation ring of $X_K(L)$ (resp. $\wtilde{\bfE}_L$) in $\wtilde{\bfE}_L$. For any $c\in\bfR_{>0}$, set $I_c\coloneqq \{x\in\bfC_p\colon v_p(x)\geq c\}$. There exists a constant $c(L/K)>0$ such that for any real numbers $c\in (0,c(L/K))$,
	\begin{thmenum}
		\item\label{it:64173} The natural morphism $\wtilde{\bfE}_L^+\lto \varprojlim_{x\mapsto x^p}\scrO_L/(I_c\cap\scrO_L)$ is an isomorphism.
		\item\label{it:65119} The image of $\bfE_{L/K}^+$ under the above isomorphism is
		$$\left\{(x_n)\in\varprojlim_{x\mapsto x^p}\scrO_L/I_c\colon x_{r_m}\in \scrO_{F_m}/(I_c\cap\scrO_{F_m})\subseteq \scrO_L/(I_c\cap \scrO_L),\ m=1,2,\cdots\right\},$$
		where $r_m\coloneqq \log_p([F_m\colon K_1])$ and $F_m$ is the $m$-th elementary extension of $L/K$.
	\end{thmenum}
\end{theorem}
If $L/K$ is a $p$-adic Lie extension, then Sen's result (cf. \cite{senRamificationPadicLie1972}) allows us to replace the filtration of elementary extensions by the filtration induced from the Lie structure of $\opn{Gal}(L/K)$, which is more convenient in practice (see, for instance, \cite[Proposition 2.1.3]{bergerSuperHolderVectorsField2024}). However, in the case where $K_\infty/K$ is not a Galois extension, it becomes necessary to work directly with the filtration by elementary extensions. The following statement is presumably well-established in the existing literature, although we have not identified a definitive reference.
\begin{proposition}\label{prop:7578}
	The image of $\bfE_{\tau,K}^+$, the valuation ring of $\bfE_{\tau,K}$, under the isomorphism in \Cref{it:64173} is $\varprojlim_{x\mapsto x^p}\scrO_{K_n}/(I_c\cap \scrO_{K_n})$, where $K_n=K\left(\pi^{1/p^n}\right)$.
\end{proposition}
\begin{proof}
	By \Cref{it:65119}, it left to verify that for any $n\in\bfN_{\geq 1}$, $K_n$ is the $n$-th elementary extension of $K_\infty/K$, which can be reduced to check that $K_{n+1}/K_n$ is $i_n$-elementary (cf. \cite[421]{caisCharacterizationStrictlyAPF2016}), where
	$$i_n\coloneqq \frac{e_K}{p-1}\cdot p^n+\frac{1}{p}$$
	and $e_K$ is the absolute ramification index of $K$.

	By \cite[Example 4.1]{caisCharacterizationStrictlyAPF2016} as well as the proof of \cite[Proposition 3.4]{caisCharacterizationStrictlyAPF2016}, this can be reduced to show that the (finite part) of the Newton polygon\footnote{Following \cite{caisCharacterizationStrictlyAPF2016}, we calculate the Newton polygon by using the natural valuation $v_K$ on $K$ (i.e. $v_K(p)=e_K$).} of the polynomial
	\begin{equation}\label{eq:12893}
		f_n(T)=\left(T+\pi^{1/p^{n+1}}\right)^p-\pi^{1/p^n}=T^p+\sum_{k=1}^{p-1}\binom{p}{k}\pi^{\frac{p-k}{p^{n+1}}}\cdot T^k+T^p
	\end{equation}
	consists of only one segment with slope $-\frac{i_n}{p^n}$. If we let $f_n(T)=\sum_{k=0}^p a_k T^k$, then \eqref{eq:12893} gives that
	$$v_K(a_k)=\begin{cases}
			\infty,                  & \text{ if }k=0,             \\
			e_K+\frac{p-k}{p^{n+1}}, & \text{ if }1\leq k\leq p-1, \\
			0,                       & \text{ if }k=p.
		\end{cases}
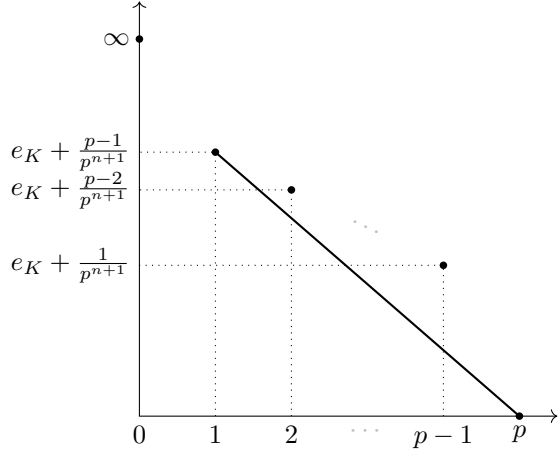
$$
	\begin{figure}[H]
		\centering
		\begin{tikzpicture}
			\draw [<->] (0,5.5) -- (0,0) -- (5.5,0);
			\node [below] at (0,0) {$0$};
			\node [below] at (1,0) {$1$};
			\node [below] at (2,0) {$2$};
			\node [below,lightgray] at (3,0) {$\cdots$};
			\node [below] at (4,0) {$p-1$};
			\node [below] at (5,0) {$p$};
			\draw [thick] (1,3.5) -- (5,0);
			\node [circle,fill,inner sep=1pt] at (1,3.5) {};
			\draw [dotted] (1,0) -- (1,3.5) --(0,3.5);
			\node [circle,fill,inner sep=1pt] at (2,3) {};
			\draw [dotted] (2,0) -- (2,3) --(0,3);
			\node [circle,fill,inner sep=1pt] at (4,2) {};
			\draw [dotted] (4,0) -- (4,2) --(0,2);
			\node [circle,fill,inner sep=1pt] at (5,0) {};
			\node [left] at (0,3.5) {$e_K+\frac{p-1}{p^{n+1}}$};
			\node [left] at (0,3) {$e_K+\frac{p-2}{p^{n+1}}$};
			\node [left] at (0,2) {$e_K+\frac{1}{p^{n+1}}$};
			\node [left] at (0,5) {$\infty$};
			\node [circle,fill,inner sep=1pt] at (0,5) {};
			\node [lightgray] at (3,2.49) {\rotatebox[origin=c]{-26.57}{$\cdots$}};
		\end{tikzpicture}
	\end{figure}
	And the result follows.
\end{proof}

\section{Super-H\"{o}lder vectors and locally analytic vectors}
\begin{definition}[cf. {\cite[Section 1]{bergerSuperHolderVectorsField2024}}]\label{def:47537}
	Let $E$ be a field of characteristic $p$, $M$ a $E$-vector space, $G$ a uniform pro-$p$ group, satisfying:
	\begin{enumerate}[label=(\roman*)]
		\item $M$ is equipped with a valuation $\opn{val}_M$ such that $\opn{val}_M(x\cdot m)=\opn{val}(m)$ holds for any $x\in E^\times$ and $m\in M$;
		\item $G$ $E$-linearly acts on $M$, and the action is isometric.
	\end{enumerate}

	Fix $e\in\bfN_{\geq 1}$, $\lambda,\mu\in\bfR$. Denote by $\calH_e^{\lambda,\mu}(G,H)$ the set of functions $f\colon G\lto M$ such that $\opn{val}_M(f(g)-f(h))\geq p^{\lambda}\cdot p^{ei}+\mu$ holds for any $i\in\bfN$ and any elements $g,h\in G$ with $g^{-1}h\in G^{(i)}$. Let
	$$\calH_e^\lambda(G,M)\coloneqq \bigcup_{\mu\in\bfR}\calH_e^{\lambda,\mu}(G,H),\ \calH_e(G,M)\coloneqq \bigcup_{\lambda\in\bfR}\calH_e^{\lambda}(G,H).$$
	Call elements in $\calH_1(G,M)$ \textbf{super-H\"{o}lder functions}. Besides that, we say $m\in M$ is a \textbf{super-H\"{o}lder vector} (with respect to the action of $G$) if the orbit map $\opn{orb}_m\colon g\longmapsto g\cdot m$ lies in $\calH_1(G,M)$. Denote by $M^{G\hyphen\opn{sh}}$ the set of super-H\"{o}lder vectors in $M$ and by $M^{G\hyphen\opn{sh},\lambda}$ (resp. $M^{G\hyphen\opn{sh},\lambda,\mu}$) the set of vectors $m\in M$ such that $\opn{orb}_m$ lies in $\calH_1^\lambda(G,M)$ (resp. $\calH_1^{\lambda,\mu}(G,M)$).
\end{definition}
\begin{lemma}[cf. {\cite[Lemma 1.2.8]{bergerSuperHolderVectorsField2024}}]
	With the above notations, let $m\in M$. Then $m\in M^{G\hyphen\opn{sh},\lambda,\mu}$ if and only if for any $i\in\bfN$ and any $g\in G^{(i)}$, one has $\opn{val}_M((g-1)\cdot m)\geq p^{\lambda}\cdot p^{ei}+\mu$.
\end{lemma}

We refer the reader to \cite[Section 1.1, Section 1.2]{bergerDecompletionCyclotomicPerfectoid2022} and \cite[Section 1.2]{bergerSuperHolderVectorsField2024} for more general properties of super-H\"{o}lder functions and vectors.
\begin{definition}\label{def:1825}
	Let $M_1$ (resp. $M_2$) be a $E$-vector space (resp. ring) as in the above definition, equipped with an isometric $E$-linear action of $\Gamma_{K}$ (resp. $\Gamma_{\rtimes,K}$). For the convenience of reader, we summarize some notations about super-H\"{o}lder vectors as follows:
	\begin{table}[H]
		\centering
		\begin{tabular}[t]{cc}
			\toprule
			Notation                                    & Definition                                                \\
			\midrule
			$M_1^{\gamma\hyphen\opn{sh}}$               & $M_1^{\Gamma_K^{(0)}\hyphen1\hyphen\opn{sh}}$             \\
			$M_1^{\gamma_k\hyphen\opn{sh},\lambda}$     & $M_1^{\Gamma_K^{(k)}\hyphen1\hyphen\opn{sh},\lambda}$     \\
			$M_1^{\gamma_k\hyphen\opn{sh},\lambda,\mu}$ & $M_1^{\Gamma_K^{(k)}\hyphen1\hyphen\opn{sh},\lambda,\mu}$ \\
			\bottomrule
		\end{tabular}
		\begin{tabular}[t]{cc}
			\toprule
			Notation                                           & Definition                                                                          \\
			\midrule
			$M_2^{\tau\hyphen\opn{sh}}$                        & $M_2^{\Gamma_{\tau,K}^{(0)}\hyphen1\hyphen\opn{sh}}$                                \\
			$M_2^{\tau_k\hyphen\opn{sh},\lambda}$              & $M_2^{\Gamma_{\tau,K}^{(k)}\hyphen1\hyphen\opn{sh},\lambda}$                        \\
			$M_2^{\tau_k\hyphen\opn{sh},\lambda,\mu}$          & $M_2^{\Gamma_{\tau,K}^{(k)}\hyphen1\hyphen\opn{sh},\lambda,\mu}$                    \\
			$M_2^{\tau\hyphen\opn{sh};\gamma=1}$               & $M_2^{\Gamma_{\tau,K}^{(0)}\hyphen1\hyphen\opn{sh}}\cap M_2^{\Gamma_K}$             \\
			$M_2^{\tau_k\hyphen\opn{sh},\lambda;\gamma=1}$     & $M_2^{\Gamma_{\tau,K}^{(k)}\hyphen1\hyphen\opn{sh},\lambda}\cap M_2^{\Gamma_K}$     \\
			$M_2^{\tau_k\hyphen\opn{sh},\lambda,\mu;\gamma=1}$ & $M_2^{\Gamma_{\tau,K}^{(k)}\hyphen1\hyphen\opn{sh},\lambda,\mu}\cap M_2^{\Gamma_K}$ \\
			\bottomrule
		\end{tabular}
	\end{table}
\end{definition}

\section{Decompletion of $\wtilde{\bfE}_{\tau,K}$}
The main goal of this section is to study the decompletion and deperfection of the perfectoid field $\wtilde{\bfE}_{\tau,K}$. Since $K_\infty/K$ is not Galois, it makes no sense to talk about the $\Gamma_{\tau,K}$-super-H\"{o}lder vectors of $\wtilde{\bfE}_{\tau,K}$. On the other hand, the work of Gao and Poyeton on overconvergent $(\varphi,\tau)$-modules (cf. \cite{gaoLOCALLYANALYTICVECTORS2021} and \cite{poyetonExtensionsLiePadiques2019}) suggests that one should consider $\Gamma_{\tau,K}$-analytic vectors in $\wtilde{\bfE}_{\rtimes,K}$ that are fixed by $\Gamma_K$. More specifically, we present the following result:
\begin{theorem}\label{thm:56067}\leavevmode
	\begin{thmenum}
		\item\label{it:59278} $\wtilde{\bfE}_{\rtimes,K}^{\tau\hyphen\opn{sh};\gamma=1}=\varphi^{-\infty}(\bfE_{\tau,K})$, where $\varphi^{-\infty}(\bfE_{\tau,K})\coloneqq \bigcup_{n\geq 0}\varphi^{-n}(\bfE_{\tau,K})$ is the perfection of $\bfE_{\tau,K}$.
		\item\label{it:21546} For any $n,k\in\bfN$, one has $$\wtilde{\bfE}_{\rtimes,K}^{\tau_k\hyphen\opn{sh},k-n+c_p;\gamma=1}=\varphi^{-n}(\bfE_{\tau,K}),$$
		where $c_p\coloneqq \log_p\frac{p}{p-1}$.
	\end{thmenum}
\end{theorem}

\subsection{Decompletion: Proof of \Cref{it:59278}}
By \cite[Proposition 1.11]{bergerDecompletionCyclotomicPerfectoid2022}, we only need to prove
\begin{equation}\label{eq:6081}
	\left(\wtilde{\bfE}_{\rtimes,K}^+\right)^{\tau\hyphen\opn{sh};\gamma=1}=\bigcup_{m\geq 0}\varphi^{-m}\left(\bfE_{\tau,K}^+\right),
\end{equation}
where $\wtilde{\bfE}_{\rtimes,K}^+\coloneqq \calO_{\widehat{K}_\rtimes}^\flat$ and $\bfE_{\tau,K}^+$ is the ring of integers of $\bfE_{\tau,K}$.

One of the inclusions in \eqref{eq:6081} is easy:
\begin{proposition}\label{prop:18986}
	For any integer $n\geq 0$, one has
	$$\varphi^{-n}\left(\bfE_{\tau,K}^+\right)\subseteq \left(\wtilde{\bfE}_{\rtimes,K}^+\right)^{\tau_k\hyphen \opn{sh},k-n+c_p;\gamma=1}\subseteq \left(\wtilde{\bfE}_{\rtimes,K}^+\right)^{\tau\hyphen \opn{sh};\gamma=1}$$
\end{proposition}
\begin{proof}
	Take any element $y=\sum_{r=0}^\infty y_r\cdot \ulpi^r$ with $y_r\in \kappa$. For any $i\in\bfN$ and any element $\tau^{m\cdot p^{k+i}}\in \Gamma_{\tau,K}^{(k+i)}$ with $m\in\bfZ_p$, one has
	$$v^\flat\left(\left(\tau^{m\cdot p^{k+i}}\right)y\right)=v^\flat\left(\sum_{r=1}^\infty \left(\varepsilon^{mr\cdot p^{k+i}}-1\right)\ulpi^r\right)\geq \min_{r\geq 1}v^\flat\left(\left(\varepsilon^{mr\cdot p^{k+i}}-1\right)\ulpi^r\right).$$
	Since \Cref{lem:61333} implies that $v^\flat\left(\varepsilon^{mr\cdot p^{k+i}}-1\right)=p^{k+i+v_p(mr)}+c_p$, one has
	$$v^\flat\left(\left(\tau^{m\cdot p^{k+i}}\right)y\right)\geq p^{k+i+c_p}+v_p(\ulpi),$$
	i.e. $y\in\left(\wtilde{\bfE}_{\rtimes,K}^+\right)^{\tau_k\hyphen\opn{sh},k+c_p,v_p(\ulpi)}$.

	Now we suppose $\varphi^n(x)\in\bfE_{\tau,K}^+\subseteq \wtilde{\bfE}_{\rtimes,K}$. Then the above calculation shows that for any $i\in\bfN$ and any $g\in \Gamma_{\tau,K}^{(k+i)}$, one has
	$$v^\flat(g\cdot x-x)=p^{-n}\cdot v^\flat (g\cdot \varphi^n(x)-\varphi^n(x))\geq p^{k-n+c_p}\cdot p^i+p^{-n}\cdot v_p(\ulpi),$$
	i.e.
	$x\in \left(\wtilde{\bfE}_{\rtimes,K}^+\right)^{\tau_k\hyphen\opn{sh},k-n+c_p}$. Since $\bigcup_{n\in\bfN}\varphi^{-n}(\bfE_{\tau,K}^+)\subseteq \wtilde{\bfE}_{\tau,K}=\wtilde{\bfE}_{\rtimes,K}^{\Gamma_K}$, the result follows.
\end{proof}
\begin{lemma}\label{lem:61333}
	For any $n\in\bfN_{\geq 1}$, we have $v^\flat\left(\varepsilon^m-1\right)=p^{v_p(m)+c_p}$, where $c_p\coloneqq \log_p\frac{p}{p-1}$.
\end{lemma}
\begin{proof}
	Note that $v^\flat\left(\varepsilon^m-1\right)=v_p\left(\lim_{k\to\infty}\left(\zeta_{p^k}^m-1\right)^{p^k}\right)$. The result follows from the fact that $\zeta_{p^k}^m$ is a $p^{k-v_p(m)}$-th primitive root of unity, and $\zeta_{p^k}^m-1$ is a uniformizer of the $p^{k-v_p(m)}$-th cyclotomic extension of $\bfQ_p$ with $p$-adic valuation $p^{v_p(m)+c_p}\cdot p^{-k}$.
\end{proof}

Compared to the cyclotomic case, the opposite inclusion in \eqref{eq:6081} is more subtle, as the Galois groups involved are more complicated, and the non-Galois nature of $K_\infty/K$ forces us to do a more low-level calculation on the image of the embedding $X_K^+(K_\infty)\hookrightarrow \wtilde{\bfE}^+$ via analysing the ramification of $K_\infty/K$ (cf. \Cref{prop:7578}), which is bypassed in the cyclotomic case via Sen's theorem. Nevertheless, we can still follow the strategy in \cite[Theorem 2.2.3]{bergerSuperHolderVectorsField2024} to prove the following proposition:
\begin{proposition}
	Let $x\in \left(\wtilde{\bfE}_{\rtimes,K}^+\right)^{\tau_k\hyphen \opn{sh},\lambda,\mu;\gamma=1}$ for some $\lambda,\mu\in\bfR$. Then there exists $n\in\bfN$ such that $\varphi^n(x)\in \bfE_{\tau,K}^+$.
\end{proposition}
\begin{proof}
	By \cite[Lemma 1.2.11]{bergerSuperHolderVectorsField2024}, one may assume that $k=0$. The condition $x\in \left(\wtilde{\bfE}_{\rtimes,K}^+\right)^{\tau_0\hyphen \opn{sh},\lambda,\mu;\gamma=1}$ implies that for any $i\in\bfN$ and any $g\in\Gamma_{\tau,K}^{(i)}$, one has
	$$v^\flat(g\cdot x-x)\geq p^\lambda\cdot p^i+\mu.$$
	Without loss of generality, we may further assume that $\lambda$ is large enough such that the inequality $p^\lambda\cdot p^i+\mu\geq 0$ holds for any $i\in\bfN$. As a result, one may take $n\in\bfN$ large enough to ensure
	\begin{equation}\label{eq:47168}
		v^\flat\left(g\cdot \varphi^n(x)-\varphi^n(x)\right)\geq p^n\cdot\left(p^\lambda\cdot p^i+\mu\right)\geq p^i
	\end{equation}
	holds for any integer $i\geq 1$.

	Since $\wtilde{\bfE}_{\rtimes,K}^{\Gamma_K}=\wtilde{\bfE}_{\tau,K}$, we can assume
	$$\varphi^n(x)=(y_k)_{k\in\bfN}\in\varprojlim_{x\mapsto x^p}\scrO_{K_\infty}/p=\wtilde{\bfE}_{\tau,K}^+\subsetneq\wtilde{\bfE}_{\rtimes,K}^+.$$
	Besides that, we set $y^{(k)}\coloneq \lim_{m\to\infty}\widehat{y}_{k+m}^{p^m}$, where $\widehat{y}_{k+m}\in\scrO_{K_\infty}$ is a lift of $y_{k+m}$. For any $i\in\bfN_{\geq 1}$ and any $g\in \Gamma_{\tau,K}^{(i)}$, one has
	$$v^\flat\left(g\cdot \varphi^n(x)-\varphi^n(x)\right)=v_p\left(\left(g\cdot \varphi^n(x)-\varphi^n(x)\right)^\sharp\right)=p^i\cdot v_p\left(\left(g\cdot \varphi^n(x)-\varphi^n(x)\right)^{(i)}\right).$$
	By \Cref{eq:47168}, this implies that $$v_p\left(\left(g\cdot \varphi^n(x)-\varphi^n(x)\right)^{(i)}\right)\geq 1.$$
	Notice that
	$$\left(g\cdot \varphi^n(x)-\varphi^n(x)\right)^{(i)}\equiv g\cdot \left(\varphi^n(x)\right)^{(i)}-\left(\varphi^n(x)\right)^{(i)}=g\cdot y^{(i)}-y^{(i)}\pmod{p},$$
	strong triangle inequality implies that
	\begin{equation}\label{eq:20215}
		v_p\left(g\cdot y^{(i)}-y^{(i)}\right)\geq 1
	\end{equation}
	holds for any $i\in\bfN_{\geq 1}$. On the other hand, since the action of $\scrG_{K_i}$ on $K_\infty$ factors through
	$$\opn{Gal}(K_\rtimes/K_i)=\opn{Gal}\left((K_i)_\rtimes/K_i\right)\cong\Gamma_{\tau,K_i}\rtimes \Gamma_{K_i},$$
	where $\Gamma_{\tau,K_i}=\Gamma_{\tau,K}^{(i)}$ by \Cref{lem:23768} and
	$$\Gamma_{K_i}=\opn{Gal}((K_i)_{p^\infty}/K_i)\xlongrightarrow[\text{\Cref{it:45121}}]{\cong}\opn{Gal}((K_i)_\rtimes/(K_i)_\infty)=\Gal(K_\rtimes/K_\infty),$$
	while the action of $\opn{Gal}(K_\rtimes/K_i)\cong \Gamma_{\tau,K}^{(i)}\rtimes \Gal(K_\rtimes/K_\infty)$ factors through its set-theoretic quotient $\Gamma_{\tau,K}^{(i)}$, one can conclude from \Cref{eq:20215} that $v_p\left(g\cdot y^{(i)}-y^{(i)}\right)\geq 1 $ holds for any $g\in\scrG_{K_i}$.

	If we fix $c\in (0,c(K_\infty/K))$ as in \Cref{thm:10788}, then by the explicit version of Ax-Sen-Tate theorem \cite[Th\'eor\`eme 1.7]{leborgneOptimisationTheoremeAxSenTate2010a}, the image of $y^{(i)}$ under the projection $\scrO_{K_\infty}\longrightarrow \scrO_{K_\infty}/(I_c\cap \scrO_{K_\infty})$ belongs to $\scrO_{K_n}/(I_c\cap \scrO_{K_n})$. As a result, the image of $\varphi^n(x)=(y_1^p,y_1,\cdots,y_i,\cdots)$ under the projection
	$$\wtilde{\bfE}_{\tau,K}=\varprojlim_{x\mapsto x^p}\widehat{\scrO}_{K_\infty}\xrightarrow[\Cref{it:64173}]{\cong}\varprojlim_{x\mapsto x^p}\scrO_{K_\infty}/(I_c\cap \scrO_{K_\infty})$$
	lies in $\varprojlim_{x\mapsto x^p}\scrO_{K_n}/(I_c\cap \scrO_{K_n})$. This finishes the proof by \Cref{prop:7578}.
\end{proof}

\subsection{Deperfection: proof of \Cref{it:21546}}
In \cite{bergerDecompletionCyclotomicPerfectoid2022}, the deperfection result in the cyclotomic case is achieved via the Tate normalized trace maps. Although we can define the normalized Tate trace maps for the Kummer extension as in the cyclotomic case, they are not ``$\Gamma_{\tau,K}$-equivariant''. On the other hand, we observe that \Cref{it:59278} is already strong enough to supersede the role of Tate normalized trace maps (i.e. \cite[Section 2.2]{bergerDecompletionCyclotomicPerfectoid2022}).

The following lemma is the adaptation of \cite[Proposition 2.4]{bergerDecompletionCyclotomicPerfectoid2022} to our Kummer setting.
\begin{lemma}\label{lem:47455}
	For any $d\in\bfR_{>0}$, one has $\left(\wtilde{\bfE}_{\rtimes,K}^+\right)^{\tau_k\hyphen \opn{sh},k+c_p+d;\gamma=1}\subset \varphi\left(\bfE_{\tau,K}^+\right)$.
\end{lemma}
\begin{proof}
	Fix $d\in\bfR_{>0}$. Take $f(T)\in\kappa\bbrac{T}$ such that $f^\prime(T)\neq 0$, i.e. $f(T)\notin \kappa\bbrac{T^p}$. We prove that $f(\ulpi)$ does not belong to $\left(\wtilde{\bfE}_{\rtimes,K}^+\right)^{\tau_k\hyphen \opn{sh},k+c_p+d;\gamma=1}$.

	For any $i\in\bfN$, let $g\coloneqq p^{k+i}\in\Gamma_{\tau,K}^{(k+i)}$. The Taylor expansion gives that
	\begin{align*}
		g\cdot f(\ulpi)= & f\left(\left(\varepsilon^{p^{k+i}}-1\right)\ulpi+\ulpi\right)                                                                              \\
		=                & f(\ulpi)+\left(\varepsilon^{p^{k+i}}-1\right)\ulpi\cdot f^\prime(\ulpi)+\left(\left(\varepsilon^{p^{k+i}}-1\right)\ulpi\right)^2\cdot a_g,
	\end{align*}
	where $a_g$ is some element in $\wtilde{\bfE}_{\rtimes,K}^+$. As a result, one can take integer $i_\mu\gg 0$ for arbitrary $\mu\in\bfR$ such that
	\begin{enumerate}
		\item $v^\flat\left(f^\prime(\ulpi)\right)<v^\flat\left(\left(\varepsilon^{p^{k+i_\mu}}-1\right)\ulpi\right)$;
		\item $v^\flat\left(f^\prime(\ulpi)\right)+v^\flat(\ulpi)+p^{k+c_p}\cdot p^{i_\mu}<p^{k+c_p+d}\cdot p^{i_\mu}+\mu$.
	\end{enumerate}
	Then the strong triangle inequality implies that
	\begin{align*}
		v^\flat\left(g\cdot f(\ulpi)-f(\ulpi)\right)= & v^\flat\left(\left(\varepsilon^{p^{k+i_\mu}}-1\right)\ulpi\cdot f^\prime(\ulpi)\right) \\
		=                                             & v^\flat\left(f^\prime(\ulpi)\right)+v^\flat(\ulpi)+p^{k+i_\mu+c_p}                     \\
		<                                             & p^{k+c_p+d}\cdot p^{i_\mu}+\mu,
	\end{align*}
	i.e. $f(\ulpi)$ does not belong to $\left(\wtilde{\bfE}_{\rtimes,K}^+\right)^{\tau_k\hyphen \opn{sh},k+c_p+d,\mu}$ for any $\mu\in\bfR$. This finishes the proof.
\end{proof}

\begin{proof}[Proof of \Cref{it:21546}]
	A mimic of \cite[Corollary 2.5]{bergerDecompletionCyclotomicPerfectoid2022}, which depends on \Cref{lem:47455}, gives us
	$$\left(\wtilde{\bfE}_{\rtimes,K}^+\right)^{\tau_k\hyphen\opn{sh},k-n+c_p}\cap \bigcup_{n\geq 0}\varphi^{-n}\left(\bfE_{\tau,K}^+\right)=\varphi^{-n}\left(\bfE_{\tau,K}^+\right).$$

	On the other hand, the decompletion result \Cref{it:59278} implies that
	$$\left(\wtilde{\bfE}_{\rtimes,K}^+\right)^{\tau_k\hyphen\opn{sh},k-n+c_p;\gamma=1}=\left(\wtilde{\bfE}_{\rtimes,K}^+\right)^{\tau_k\hyphen\opn{sh},k-n+c_p}\cap \bigcup_{n\geq 0}\varphi^{-n}\left(\bfE_{\tau,K}^+\right).$$
	The result follows.
\end{proof}

\section{Descent of $(\varphi,\tau)$-modules over perfectoid period ring}
In this section, we generalize the decompletion and deperfection of $(\varphi,\Gamma)$-modules by Berger and Rozensztajn (cf. \cite[Corollary 3.11]{bergerDecompletionCyclotomicPerfectoid2022}) to the setting of $(\varphi,\tau)$-modules.
\begin{theorem}\label{thm:7812}
	Let $M$ be an \'etale $(\varphi,\tau)$-module over $(\bfE_{\tau,K},\wtilde{\bfE}_{\rtimes,K})$ and $k\in\bfN$. Then
	$$\left(\wtilde{\bfE}_{\rtimes,K}\otimes_{\bfE_{\tau,K}}M\right)^{\tau\hyphen\opn{sh};\gamma=1}=\varphi^{-\infty}(\bfE_{\tau,K})\otimes_{\bfE_{\tau,K}}M.$$
	More precisely, for any $n\in\bfN$, one has
	\begin{equation}\label{eq:5230}
		\left(\wtilde{\bfE}_{\rtimes,K}\otimes_{\bfE_{\tau,K}}M\right)^{\tau_k\hyphen\opn{sh},k-n+c_p;\gamma=1}=\varphi^n(\bfE_{\tau,K})\otimes_{\bfE_{\tau,K}}M.
	\end{equation}
\end{theorem}

In the rest of this section, we denote by $M$ an \'etale $(\varphi,\tau)$-module over $(\bfE_{\tau,K},\wtilde{\bfE}_{\rtimes,K})$ of rank $d$ and fix $k\in\bfN$.

\subsection{Selection of lattices and induced valuations}\label{sec:37644}
Although there should be no fundamental distinction between the descent theory of $(\varphi,\tau)$-modules and that of $(\varphi,\Gamma)$-modules, the non-Galois nature of Kummer extension $K_\infty/K$ brings in some technical obstacles.

For instance, \cite[Corollary 3.11]{bergerDecompletionCyclotomicPerfectoid2022} depends on the existence of a $\Gamma_K$-stable $\mathbf{E}_K^+$-lattice within any $\mathbf{E}_K$-vector space $D$ equipped with a semi-linear $\Gamma_K$-action, which will induce a valuation of $D$ such that $\Gamma_K$ acts isometrically on $D$. In contrast, for general $\mathbf{E}_{\tau,K}$-vector space $M$, there is no well-defined $\Gamma_{\tau,K}$-action on it, let alone the $\Gamma_{\tau,K}$-stable $\mathbf{E}_{\tau,K}^+$-lattice in $M$. To overcome this difficulty, we consider two different valuations on $\wtilde{\bfE}_{\rtimes,K}\otimes_{\bfE_{\tau,K}}M$:
\begin{enumerate}
	\item Fix a $\bfE_{\tau,K}$-base $\Phi_\tau=(m_1,\cdots,m_d)$ of $M$. It induces a valuation $v_\tau$ on $\wtilde{\bfE}_{\rtimes,K}\otimes_{\bfE_{\tau,K}}M$:
	      $$v_\tau\colon \sum_{j=1}^d a_j\cdot m_j\longmapsto \min_{1\leq j\leq d}v^\flat(a_j).$$
	      In other words, this is the valuation that coming from the underlying $\varphi$-module, without considering the $\Gamma_{\tau,K}$-action.
	\item Fix a $\Gamma_{\tau,K}$-stable $\wtilde{\bfE}_{\rtimes,K}^+$-lattice $\scrL=\bigoplus_{j=1}^d\wtilde{\bfE}_{\rtimes,K}^+\cdot w_j$ of $\wtilde{\bfE}_{\rtimes,K}\otimes_{\bfE_{\tau,K}}M$ (cf. \Cref{lem:50238}), which will induce a valuation $\wtilde{v}$ on $\wtilde{\bfE}_{\rtimes,K}\otimes_{\bfE_{\tau,K}}M$ such that $\Gamma_{\tau,K}$ acts isometrically on it (\Cref{lem:19133}). This valuation allows us to consider super-H\"{o}lder vectors in $\wtilde{\bfE}_{\rtimes,K}\otimes_{\bfE_{\tau,K}}M$, at the cost that it might not coming from the underlying $\varphi$-module.
\end{enumerate}
It is well known that these two valuations are equivalent (cf. \cite[Theorem 1.3.6]{kedlayaPadicDifferentialEquations2010}), i.e. there exists some constant $C>0$ satisfying $|v_\tau(x)-\wtilde{v}(x)|\leq C$ for any $x\in \wtilde{\bfE}_{\rtimes,K}\otimes_{\bfE_{\tau,K}}M$. We use these two valuations alternatively to get around the non-Galois difficulty.

\begin{lemma}\label{lem:50238}
	Let $M$ be an \'etale $(\varphi,\tau)$-module over $(\bfE_{\tau,K},\wtilde{\bfE}_{\rtimes,K})$. There exists a $\Gamma_{\tau,K}$-stable $\wtilde{\bfE}_{\rtimes,K}^+$-lattice $\scrL$ in $\wtilde{\bfE}_{\rtimes,K}\otimes_{\bfE_{\tau,K}}M$.
\end{lemma}
\begin{proof}
	If we write $M=\bfD_\tau(V)$ for some $\bfF_p$-representation $V$ of $\scrG_K$, then it is well-known that one has $\scrG_K$-equalvariant isomorphisms
	\begin{equation}\label{eq:569}
		M=\wtilde{\bfE}_{\rtimes,K}\otimes_{\bfE_{\tau,K}}\left(\bfE_{\tau,K}^{\opn{sep}}\otimes_{\bfF_p}V\right)^{\scrG_{K_\infty}}\cong \left(\wtilde{\bfE}\otimes_{\bfF_p}V\right)^{\scrG_{K_\rtimes}}\cong \tilde{\bfE}_{\rtimes,K}\otimes_{\bfE_{\rtimes,K}}\bfD_\rtimes(V),
	\end{equation}
	where $\bfD_\rtimes(V)\coloneqq\left(\bfE_{\rtimes,K}^{\opn{sep}}\otimes_{\bfF_p}V\right)^{\scrG_{K_\rtimes}}$ is the analogue of \'etale $(\varphi,\Gamma)$-module for $K_\rtimes/K$ over $\bfE_{\rtimes,K}$ that associated to $V$. By applying \cite[Proposition 3.4]{bergerDecompletionCyclotomicPerfectoid2022} to $\bfD_{\rtimes(V)}$, we get a $\bfE_{\rtimes,K}^+$-lattice in $\bfD_{\rtimes}(V)$ that stable under the action of $\Gamma_{\rtimes,K}$. It induces the expected lattice through the isomorphisms in \Cref{eq:569}.
\end{proof}
\begin{definition}
	For any matrix $Q=(q_{ij})_{1\leq i,j \leq d}$ in $\opn{M}_d(\wtilde{\bfE})$, we define the valuation of $Q$ as $v^\flat(Q)\coloneqq \min_{1\leq i,j\leq d}v^\flat(q_{ij})$.
\end{definition}
\begin{lemma}\label{lem:19133}
	The action of $\Gamma_{\tau,K}$ on $\wtilde{\bfE}_{\rtimes,K}\otimes_{\bfE_{\tau,K}}M$ is isometric with respect to the valuation $\wtilde{v}$.
\end{lemma}
\begin{proof}
	The $\Gamma_{\tau,K}$-stability of the lattice $\scrL$ implies that for any $g\in \Gamma_{\tau,K}$, the matrix of $g$ with respect to the base of $\scrL$, which we denote by $\wtilde{\opn{Mat}}(Q)$, satisfies that $v^\flat(\wtilde{\opn{Mat}}(Q))\geq 0$. As a result, for any $x=\sum_{i=1}^d x_i\cdot w_i\in \wtilde{\bfE}_{\rtimes,K}\otimes_{\bfE_{\tau,K}}M$ and any $g\in\Gamma_{\tau,K}$, one has
	\begin{align*}
		\wtilde{v}(g\cdot x)= & \wtilde{v}\left(\sum_{i=1}^d\left(\sum_{j=1}^d g(a_j)\wtilde{\opn{Mat}}(g)_{ij}\right)\cdot w_i\right)              \\
		=                     & \min_{1\leq i\leq d}v^\flat\left(\sum_{j=1}^d g(a_j)\wtilde{\opn{Mat}}(g)_{ij}\right)                               \\
		\geq                  & \min_{1\leq i\leq d}\min_{1\leq j\leq d}\left(v^\flat(g(a_j))+v^\flat\left(\wtilde{\opn{Mat}}(g)_{ij}\right)\right) \\
		=                     & \min_{1\leq j\leq d}v^\flat(a_j)+v^\flat\left(\wtilde{\opn{Mat}}(g)_{ij}\right)                                     \\
		\geq                  & \wtilde{v}(x).
	\end{align*}
	By writting $x=g^{-1}\cdot (g\cdot x)$ and applying the above inequality to $g^{-1}$, one gets the equality $\wtilde{v}(g\cdot x)=\wtilde{v}(x)$.
\end{proof}

\subsection{Calculation of super-H\"{o}lder vectors}
In this section, building upon the framework established in \cite[Section 3.2]{bergerDecompletionCyclotomicPerfectoid2022}, we determine the set of super-Hölder vectors
$$\left(\wtilde{\bfE}_{\rtimes,K}\otimes_{\bfE_{\tau,K}}M\right)^{\tau_k\hyphen\opn{sh},k-n+c_p;\gamma=1}$$
for every $n\in\bfN$. Notations introduced in \Cref{sec:37644} will be freely used without further mention, and we denote by $\opn{Mat}(g)$ the matrix of $g\in \Gamma_{\tau,K}$ with respect to the base $\Phi_\tau$ of $M$.

The following lemma allows us to reduce the calculation of super-H\"{o}lder vectors to a set of base elements of $M$:
\begin{lemma}[cf. {\cite[Lemma 3.5]{bergerDecompletionCyclotomicPerfectoid2022}\textsuperscript{\textdagger}}]\label{lem:17409}%
	The set $\left(\wtilde{\bfE}_{\rtimes,K}\otimes_{\bfE_{\tau,K}}M\right)^{\tau_k\hyphen\opn{sh},\lambda+c_p;\gamma=1}$ is a $\bfE_{\tau,K}$-vector subspace of $\wtilde{\bfE}_{\rtimes,K}\otimes_{\bfE_{\tau,K}}M$ when $\lambda\leq k$. In particular, the set $\left(\wtilde{\bfE}_{\rtimes,K}\otimes_{\bfE_{\tau,K}}M\right)^{\tau\hyphen\opn{sh};\gamma=1}$ is a $\bfE_{\tau,K}$-vector subspace of $\wtilde{\bfE}_{\rtimes,K}\otimes_{\bfE_{\tau,K}}M$.

\end{lemma}
As observed by Berger and Rozensztajn in \cite[Proposition 3.2]{bergerDecompletionCyclotomicPerfectoid2022}, the decompletion (as well as deperfection) of a $(\varphi,\Gamma)$-modules $D$ over $\wtilde{\bfE}_K$ can be deduced from showing $D$ itself is super-H\"{o}lder. In our $(\varphi,\tau)$-modules setting, we have the following analogue:
\begin{proposition}
	\Cref{thm:7812} is implied by the condition that $M\subseteq \left(\wtilde{\bfE}_{\rtimes,K}\otimes_{\bfE_{\tau,K}M}\right)^{\tau_k\hyphen\opn{sh},k+c_p}$.
\end{proposition}
\begin{proof}
	Fix $n\in\bfN$. Similar to \Cref{lem:17409}, we can show that $\varphi^{-n}\left(\bfE_{\tau,K}\right)\otimes_{\bfE_{\tau,K}}M$ is a $\varphi^{-n}\left(\bfE_{\tau,K}\right)$-vector subspace of $\left(\wtilde{\bfE}_{\rtimes,K}\otimes_{\bfE_{\tau,K}}M\right)^{\tau_k\hyphen\opn{sh},k-n+c_p;\gamma=1}$. Since $\left(\wtilde{\bfE}_{\rtimes,K}\otimes_{\bfE_{\tau,K}}M\right)^{\tau_k\hyphen\opn{sh},k-n+c_p;\gamma=1}$ contains $M$ by assumption, it also contains $\varphi^{-n}\left(\bfE_{\tau,K}\right)\otimes_{\bfE_{\tau,K}}M$.

	To show the opposite inclusion, we take
	$$x\coloneqq \sum_{j=1}^d x_j\cdot m_j\in \left(\wtilde{\bfE}_{\rtimes,K}\otimes_{\bfE_{\tau,K}}M\right)^{\tau_k\hyphen\opn{sh},k-n+c_p;\gamma=1},\ x_1,\cdots,x_d\in\wtilde{\bfE}_{\rtimes,K},$$
	and write $g\cdot x=\sum_{j=1}^d f_j(g)\cdot m_j$ for any $g\in \Gamma_{\tau,K}^{(k)}$. Then one has $x_j\in \wtilde{\bfE}_{\rtimes,K}^{\Gamma_K}$ for $j=1,\cdots,d$. Suppose $m_j\in \left(\wtilde{\bfE}_{\rtimes,K}\otimes_{\bfE_{\tau,K}}M\right)^{\tau_k\hyphen\opn{sh},k-n+c_p,\mu_j}$ for certain $\mu_j\in\bfR$, $j=1,\cdots,d$. Then for any $i\in\bfN$, $j=1,\cdots,d$ and $g,h\in\Gamma_{\tau,K}^{(k+i)}$ satisfying $g\cdot h^{-1}\in\Gamma_{\tau,K}^{(k+i)}$, one has
	\begin{align*}
		v^\flat\left(f_j(g)-f_j(h)\right)\geq & v_\tau(g\cdot x-h\cdot x)       \\
		\geq                                  & \wtilde{v}(g\cdot x-h\cdot x)-C \\
		\geq                                  & p^{k-n+c_p}\cdot p^i+\mu-C,
	\end{align*}
	i.e. the function
	$$g\longmapsto \begin{pmatrix}
			f_1(g) & 0 & \cdots & 0 \\
			f_2(g) & 0 & \cdots & 0 \\
			\vdots &   &        &   \\
			f_d(g) & 0 & \cdots & 0
		\end{pmatrix}$$
	belongs to $\calH_1^{k-n+c_p}\left(\Gamma_{\tau,K}^{(k)},\opn{M}_s\left(\wtilde{\bfE}_{\rtimes,K}\right)\right)$. On the other hand, by the inclusion
	$$M\subseteq\left(\wtilde{\bfE}_{\rtimes,K}\otimes_{\bfE_{\tau,K}}M\right)^{\tau_k\hyphen\opn{sh},k+c_p}\subseteq \left(\wtilde{\bfE}_{\rtimes,K}\otimes_{\bfE_{\tau,K}}M\right)^{\tau_k\hyphen\opn{sh},k-n+c_p}$$
	and \Cref{prop:15428}, the function $g\longmapsto \opn{Mat}(g)$ also belongs to $\calH_1^{k-n+c_p}\left(\Gamma_{\tau,K},\opn{M}_n\left(\wtilde{\bfE}_{\rtimes,K}\right)\right)$.

	Since $\Gamma_{\tau,K}^{(k)}\cong\bfZ_p$ is compact and $\opn{det}\opn{Mat}(g)\neq 0$ holds for any $g\in \Gamma_{\tau,K}^{(k)}$, the image of the continuous map $g\longmapsto \opn{det}\opn{Mat}(g)$ is compact and does not contain 0. This allows us to apply \cite[Proposition 1.4 (1) (4)]{bergerDecompletionCyclotomicPerfectoid2022} to show that the map
	$$g\longmapsto \begin{pmatrix}
			g(x_1) & 0 & \cdots & 0 \\
			g(x_2) & 0 & \cdots & 0 \\
			\vdots &   &        &   \\
			g(x_d) & 0 & \cdots & 0
		\end{pmatrix}=\opn{Mat}(g)^{-1}\cdot \begin{pmatrix}
			f(x_1) & 0 & \cdots & 0 \\
			f(x_2) & 0 & \cdots & 0 \\
			\vdots &   &        &   \\
			f(x_d) & 0 & \cdots & 0
		\end{pmatrix}$$
	belongs to $\calH_1^{k-n+c_p}\left(\Gamma_{\tau,K},\opn{M}_n\left(\wtilde{\bfE}_{\rtimes,K}\right)\right)$, and consequently $\opn{orb}_{x_j}\in \calH_1^{k-n+c_p}\left(\Gamma_{\tau,K},\wtilde{\bfE}_{\rtimes,K}\right)$ for $j=1,\cdots,d$. This indicates that $x_j$ belongs to $\left(\wtilde{\bfE}_{\rtimes,K}\right)^{\tau_k\hyphen\opn{sh},k-n+c_p;\gamma=1}$, which equals to $\varphi^{-n}\left(\bfE_{\tau,K}\right)$ by \Cref{it:21546}.
\end{proof}
\begin{proposition}[cf. {\cite[Proposition 3.6]{bergerDecompletionCyclotomicPerfectoid2022}\textsuperscript{\textdagger}}]\label{prop:15428}
	Let $\lambda\in\bfR$ with $\lambda\leq k$. Then $M\subseteq\left(\wtilde{\bfE}_{\rtimes,K}\otimes_{\bfE_{\tau,K}}M\right)^{\tau_k\hyphen\opn{sh},\lambda+c_p;\gamma=1}$ if and only if the map
	$$\Gamma_{\tau,K}^{(k)}\lto \opn{M}_n\left(\wtilde{\bfE}_{\rtimes,K}^+\right),\ g\longmapsto \opn{Mat}(g)$$
	belongs to $\calH_1^{\lambda+c_p}\left(\Gamma_{\tau,K}^{(k)},\opn{M}_n\left(\wtilde{\bfE}_{\rtimes,K}\right)\right)$.
\end{proposition}
\begin{proof}
	By \Cref{lem:17409}, the condition $M\subseteq\left(\wtilde{\bfE}_{\rtimes,K}\otimes_{\bfE_{\tau,K}}M\right)^{\tau_k\hyphen\opn{sh},\lambda+c_p;\gamma=1}$ is equivalent to that for any $m_j$ in the base $\Phi_\tau$ of $M$, there exists $\mu_j\in\bfR$ such that
	\begin{equation}\label{eq:23490}
		m_j\in \left(\wtilde{\bfE}_{\rtimes,K}\otimes_{\bfE_{\tau,K}}M\right)^{\tau_k\hyphen\opn{sh},\lambda+c_p,\mu_j}.
	\end{equation}
	If \eqref{eq:23490} holds, then for any $i\in\bfN$ and any $g,h\in \Gamma_{\tau,K}^{(k)}$ such that $g\cdot h^{-1}\in\Gamma_{\tau,K}^{(k+i)}$, one has
	\begin{align*}
		v^\flat\left(\opn{Mat}(g)-\opn{Mat}(h)\right)= & \min_{1\leq j\leq d}\left(\min_{1\leq i \leq d} \left(v^\flat\left(\opn{Mat}(g)_{ij}-\opn{Mat}(h)_{ij}\right)\right)\right) \\
		=                                              & \min_{1\leq j\leq d}v_\tau\left(g\cdot m_j-h\cdot m_j\right)                                                                \\
		\geq                                           & \min_{1\leq j\leq d}\wtilde{v}\left(g\cdot m_j-h\cdot m_j\right)-C                                                          \\
		\geq                                           & p^{\lambda+c_p}\cdot p^{k+i}+\min_{1\leq j\leq d}\mu_j-C,
	\end{align*}
	which implies that the function $g\longmapsto \opn{Mat}(g)$ belongs to $\calH_1^{\lambda+c_p}\left(\Gamma_{\tau,K}^{(k)},\opn{M}_n\left(\wtilde{\bfE}_{\rtimes,K}\right)\right)$. The opposite direction follows from a similar argument.
\end{proof}
Finally, we complete the proof of \Cref{thm:7812} by showing $M$ is super-H\"{o}lder, which is the analogue of \cite[Proposition 3.9]{bergerDecompletionCyclotomicPerfectoid2022} in our setting.
\begin{proposition}
	One has $M\subseteq \left(\wtilde{\bfE}_{\rtimes,K}\otimes_{\bfE_{\tau,K}}M\right)^{\tau_k\hyphen\opn{sh},k+c_p}$.
\end{proposition}
\begin{proof}
	Since super-Hölder vectors of the $(\varphi,\tau)$-module \( M \) makes sense only upon extending its scalars to the field \(\widetilde{\mathbf{E}}_{\rtimes,K}\), it is essential to exercise caution regarding the specific fields—namely, \(\mathbf{E}_{\tau,K}\) or \(\widetilde{\mathbf{E}}_{\rtimes,K}\)—over which the coefficients of the matrices appearing in the proof are defined.

	For any $s\in\bfN$, denote by $\opn{Mat}^s(\varphi)$ the matrix of the Frobenius $\varphi$ on $\wtilde{\bfE}_{\rtimes,K}\otimes_{\bfE_{\tau,K}}M$ with respect to the base $\Phi_\tau^s\coloneqq\left\{\ulpi^s\cdot m_1,\cdots,\ulpi^s\cdot m_d\right\}\subseteq M$. Since $M$ is stable under $\varphi$, $\opn{Mat}^s(\varphi)$ lies in $\opn{M}_n\left(\bfE_{\tau,K}\right)$. On the other hand, the identity
	$$\opn{Mat}(\varphi)\cdot \varphi\left(\ulpi^sI_n\right)=\ulpi^sI_n\cdot \opn{Mat}^s(\varphi)$$
	allows us to take integer $s\gg 0$ such that $\opn{Mat}^s(\varphi) \in \opn{M}_n\left(\bfE_{\tau,K}^+\right)$. By abuse of notations, we replace $\Phi_\tau$ by $\Phi_\tau^s$ and view $P\coloneqq \opn{Mat}(\varphi)$ as an element in $\opn{M}_n\left(\bfE_{\tau,K}^+\right)$. Take $r\in\bfN_{\geq 1}$ such that $\ulpi^rP^{-1}\in\ulpi\opn{M}_n\left(\bfE_{\tau,K}^+\right)$. By continuity of the map $\Gamma_{\tau,K}\lto \opn{GL}_d\left(\wtilde{\bfE}_{\rtimes,K}\right)$, there exists $l\in\bfN_{\geq k}$ such that $v^\flat\left(\opn{Mat}(g)-\opn{Id}\right)\geq r\cdot v^\flat(\ulpi)$. Write $\opn{Mat}(g)=\opn{Id}+\ulpi^rH_g$ with $H_g\in\opn{M}_d\left(\wtilde{\bfE}_{\rtimes,K}^+\right)$.

	Write $Q_g\coloneqq \ulpi^{r(p-1)}\left(g\cdot P\right)^{-1}$. It is easy to see that $Q_g\in \ulpi\opn{M}_d\left(\wtilde{\bfE}_{\rtimes,K}^+\right)$. The commutation relation between $\varphi$ and $\Gamma_{\tau,K}$ gives $P\varphi(\opn{Mat}(g)) = \opn{Mat}(g)(g\cdot P)$ for every $g\in\Gamma_{\tau,K}^{(l)}$. As a result, one has
	$$P(g\cdot P)^{-1}-\opn{Id}=\ulpi^r\left(H_g-P\varphi(H_g)Q_g\right),$$
	implying $P(g\cdot P)^{-1}-\opn{Id}\in \ulpi^r\opn{M}_d\left(\wtilde{\bfE}_{\rtimes,K}^+\right)$. If we set
	\begin{align*}
		f_j\colon \Gamma_{\tau,K}^{(l)} & \lto\opn{M}_d\left(\wtilde{\bfE}_{\rtimes,K}^+\right),                                                                              \\
		g                               & \longmapsto  \begin{cases}
			                                               H_g-P\varphi(H_g)Q_g=\ulpi^{-r}\left(P(g\cdot P)^{-1}-\opn{Id}\right),                           & \text{if } j=0;     \\
			                                               P\varphi(P)\cdots\varphi^{j-1}(P)\varphi^j(f_0(g))\cdot \varphi^{j-1}(Q_g)\cdots\varphi(Q_g)Q_g, & \text{if } j\geq 1,
		                                               \end{cases}
	\end{align*}
	then \cite[Proposition 1.11]{bergerDecompletionCyclotomicPerfectoid2022}, \Cref{prop:18986} and \cite[Lemma 3.10]{bergerDecompletionCyclotomicPerfectoid2022}\textsuperscript{\textdagger} indicate that there exists $\mu\in\bfR$ such that functions $g\longmapsto Q_g$ and $\{f_j\}_{j\geq 0}$ all belong to $\calH_1^{l+c_p,\mu}\left(\Gamma_{\tau,K}^{(l)},\opn{M}_d\left(\wtilde{\bfE}_{\rtimes,K}^+\right)\right)$, and the summation $T_f\coloneqq \sum_{j=0}^\infty f_j$ converges in $\calH_1^{l+c_p,\mu}\left(\Gamma_{\tau,K}^{(l)},\opn{M}_d\left(\wtilde{\bfE}_{\rtimes,K}^+\right)\right)$. Since $T_f(g)=H_g$, this implies that the map $g\longmapsto H_g$, and consequently $g\longmapsto \opn{Mat}(g)=\opn{Id}+\ulpi^rH_g$, belongs to $\calH_1^{l+c_p,\mu}\left(\Gamma_{\tau,K}^{(l)},\opn{M}_d\left(\wtilde{\bfE}_{\rtimes,K}^+\right)\right)$. Then \Cref{prop:15428} shows that $M\subseteq \left(\wtilde{\bfE}_{\rtimes,K}\otimes_{\bfE_{\tau,K}}M\right)^{\tau_l,l+c_p;\gamma=1}$. The result follows from \cite[Lemma 1.2.11]{bergerSuperHolderVectorsField2024}.

\end{proof}
\backmatter
\printbibliography
\end{document}